\documentclass[12pt,leqno,oneside]{amsproc}
\usepackage[T2A]{fontenc}
\usepackage[utf8]{inputenc}
\usepackage[russian]{babel}
\usepackage{amssymb}
\allowdisplaybreaks

\renewcommand{\subsection}{\refstepcounter{subsection}%
\par\bigskip\noindent\textbf{\upshape\arabic{subsection}. }}
\renewcommand{\subsubsection}{\refstepcounter{subsubsection}%
\par\medskip\noindent\textbf{\upshape\arabic{subsection}.%
\arabic{subsubsection}.\ }}
\renewcommand{\paragraph}{\refstepcounter{paragraph}%
\par\smallskip\noindent\textbf{\upshape\arabic{subsection}.%
\arabic{subsubsection}.\arabic{paragraph}.\ }}

\numberwithin{equation}{subsection}

\makeatletter
\renewcommand{\@makefnmark}{\hbox{\small\mathsurround=0cm%
${}\hspace{0.04cm}{}^{\@thefnmark)}$}}
\renewcommand{\@makefntext}[1]{\parindent=1em\noindent\hbox to 1.8em{%
\hss${}^{\@thefnmark)}$}\,#1}
\makeatother

\title{О сепарабельных направленностях в конструктивных топологических
пространствах}
\author{А.~А.~Владимиров\footnote{Работа поддержана РФФИ, грант \No~09-06-00125.}}

\begin{document}
\renewcommand{\proofname}{\upshape Д\:о\:к\:а\:з\:а\:т\:е\:л\:ь\:с\:т\:в\:о}
\pagestyle{plain}
\begin{abstract}
В заметке устанавливаются некоторые достаточные условия, при которых сходимость
направленности точек конструктивного топологического пространства может быть
сведена к сходимости всех её регулярных подпоследовательностей. В качестве
частного случая выводится известный результат об эквивалентности сильной
и слабой конструктивной интегрируемости по Риману.
\end{abstract}
\begin{flushleft}
УДК~510.25
\end{flushleft}
\maketitle

\section{Введение}\label{pt:1}
\subsection
Из "`классического"' математического анализа хорошо известно \cite[761]{Fih}
следующее утверждение:

\subsubsection\label{prop:1:1}
{\itshape Пусть направленное множество \(\mathfrak D\) обладает хотя бы одной
конфинальной последовательностью. Тогда сходимость произвольно фиксированной
числовой направленности \(\{f_{\alpha}\}_{\alpha\in\mathfrak D}\) равносильна
сходимости всех последовательностей вида \(\{f_{x_n}\}_{n\in\mathbb N}\),
где \(\{x_n\}_{n\in\mathbb N}\) "--- произвольная конфинальная последовательность
точек направленного множества \(\mathfrak D\).
}

\medskip
К числу простейших частных случаев этого утверждения относится хрестоматийная
теорема об эквивалентности двух определений предела числовой функции
\cite[53]{Fih}.

В конструктивном математическом анализе утверждение~\ref{prop:1:1},
рассматриваемое в полной его общности, теряет силу. В частности, хорошо известны
(см., например, \cite[Теорема~4.4]{Orev:1967}) примеры \(\text{п}\)-непрерывных
отображений несепарабельных конструктивных метрических пространств, не являющихся
непрерывными. Тем не менее, некоторые частные случаи утверждения~\ref{prop:1:1}
сохраняют своё значение и в рамках конструктивного направления. Так, основной
результат работы \cite[Гл.~3]{Ku:1970} состоит, по существу, в установлении
конструктивной верности утверждения~\ref{prop:1:1} применительно к направленности
одномерных интегральных сумм Римана.

Целью настоящей заметки является установление конструктивной верности
утверждения~\ref{prop:1:1} для случая \emph{сепарабельных} направленностей точек
топологических пространств некоторого достаточно широкого класса. Точные
определения соответствующих понятий будут даны в дальнейшем.

\subsection
Все формулируемые ниже суждения и логические переходы будут пониматься нами
с точки зрения ступенчатой семантической системы А.~А.~Маркова (см., например,
\cite{VD:2009}). Термин "<множество"> будет рассматриваться в качестве сокращения
для термина "<формула языка \(\text{\upshape Я}_{\omega+1}\), не имеющая отличных
от \(|{\lozenge}\) параметров">. Под \emph{порождаемыми} множествами будут
при этом пониматься множества, заданные средствами языка \(\text{\upshape Я}_1\)
\cite[\S~12]{VD:2009}.

\section{Основные результаты}\label{pt:2}
\subsection\label{pt:2:1}
Под \emph{топологическим пространством} мы в дальнейшем будем понимать тройку
\(\mathfrak T\rightleftharpoons\langle\mathfrak S;\,\mathfrak Q;\,
\mathfrak I\rangle\) из множества точек \(\mathfrak S\), множества базисных
окрестностей \(\mathfrak Q\) и установленного между точками и базисными
окрестностями отношения инцидентности \(\mathfrak I\), удовлетворяющего
стандартным топологическим аксиомам (см., например, \cite{FDD:1970},
\cite{Cher:1972}). При этом мы, следуя обычной традиции словоупотребления,
не будем проводить различия между базисной окрестностью и множеством инцидентных
этой окрестности точек. Аналогичным образом, точки носителя \(\mathfrak S\)
топологического пространства \(\mathfrak T\) мы в дальнейшем будем сокращённо
называть "<точками пространства \(\mathfrak T\)">.

Понятие \emph{регулярной пары} базисных окрестностей мы вводим следующим образом:

\subsubsection\label{prop:2:1:1}
{\itshape Пара \(\langle U;\,V\rangle\) базисных окрестностей топологического
пространства \(\mathfrak T\) называется регулярной, если любая точка \(x\in
\mathfrak T\) или принадлежит окрестности \(U\), или не является точкой
прикосновения окрестности \(V\).
}

\subsection
Имеет место следующий простой факт:

\subsubsection\label{prop:2:2:1}
{\itshape Пусть последовательность \(\{\xi_n\}_{n\in\mathbb N}\) точек
топологического пространства \(\mathfrak T\) обладает пределом
\(\xi\in\mathfrak T\), и пусть \(\langle U;\,V\rangle\) "--- регулярная пара
базисных окрестностей точки \(\xi\). Тогда или при любом выборе номера
\(n\in\mathbb N\) точка \(\xi_n\) принадлежит окрестности \(U\), или осуществим
номер \(m\in\mathbb N\), для которого точка \(\xi_m\) не является точкой
прикосновения окрестности \(V\).
}

\begin{proof}
Согласно определению предела, осуществим номер \(N\in\mathbb N\), начиная с которого
все точки рассматриваемой последовательности лежат в окрестности \(V\).
Согласно определению \ref{prop:2:1:1}, или для всех номеров \(n\leqslant N\) точки
\(\xi_n\) принадлежат окрестности \(U\), или хотя бы для одного номера
\(m\leqslant N\) точка \(\xi_m\) не является точкой прикосновения окрестности \(V\).
В первом случае при любом выборе номера \(n\in\mathbb N\) точка \(\xi_n\)
принадлежит окрестности \(U\). Во втором случае осуществим номер \(m\in\mathbb N\),
для которого точка \(\xi_m\) не является точкой прикосновения окрестности \(V\).
\end{proof}

\subsection
Понятия \emph{конфинальной} и \emph{регулярной} последовательности точек
направленного множества мы вводим следующим образом:

\subsubsection\label{prop:2:3:1}
{\itshape Последовательность \(\{z_n\}_{n\in\mathbb N}\) точек направленного
множества \(\mathfrak D\) называется конфинальной, если для любой точки
\(\alpha\in\mathfrak D\) найдётся номер \(N\in\mathbb N\) со свойством
\[
	(\forall n\geqslant N)\qquad \alpha\leqslant z_n.
\]
}

\subsubsection
{\itshape Последовательность \(\{x_n\}_{n\in\mathbb N}\) точек направленного
множества \(\mathfrak D\) называется регулярной относительно конфинальной
последовательности \(\{z_n\}_{n\in\mathbb N}\), если при любом выборе номера
\(n\in\mathbb N\) выполняется соотношение \(z_n\leqslant x_n\).
}

\subsection
Понятие \emph{сепарабельной} направленности точек топологического пространства
мы вводим следующим образом:

\subsubsection\label{prop:2:2:2}
{\itshape Направленность \(\{f_{\alpha}\}_{\alpha\in\mathfrak D}\) точек
топологического пространства \(\mathfrak T\) называется сепарабельной,
если при любом выборе индекса \(\alpha\in\mathfrak D\) осуществимы
индекс \(\beta\geqslant\alpha\) и порождаемое множество \(\mathfrak F_{\alpha}
\subseteq\mathfrak D\), удовлетворяющие следующим условиям:
\paragraph При любом выборе индекса \(\gamma\in\mathfrak F_{\alpha}\) выполняется
неравенство \(\alpha\leqslant\gamma\).
\paragraph При любом выборе индекса \(\beta'\geqslant\beta\) точка \(f_{\beta'}\)
является точкой прикосновения множества \(\{\xi\in\mathfrak T\;:\;
(\exists\gamma\in\mathfrak F_{\alpha})\quad \xi\eqcirc f_{\gamma}\}\).
}

\subsection
Основной результат настоящей заметки состоит в следующем:

\subsubsection\label{cxef:prop}
{\itshape Пусть \(\mathfrak T\) "--- топологическое пространство, в котором
для всяких точки \(\xi\in\mathfrak T\) и её базисной окрестности \(U\)
осуществима базисная окрестность \(V\) той же точки, образующая вместе с \(U\)
регулярную пару \(\langle U;\,V\rangle\). Пусть направленность
\(\{f_{\alpha}\}_{\alpha\in\mathfrak D}\) точек пространства \(\mathfrak T\)
является сепарабельной. Пусть также неубывающая конфинальная последовательность
\(\{z_n\}_{n\in\mathbb N}\) точек направленного множества \(\mathfrak D\) такова,
что для любой регулярной относительно неё последовательности \(\{x_n\}_{%
n\in\mathbb N}\) соответствующая последовательность \(\{f_{x_n}\}_{n\in\mathbb N}\)
точек пространства \(\mathfrak T\) является сходящейся. Тогда направленность
\(\{f_{\alpha}\}_{\alpha\in\mathfrak D}\) имеет предел.
}

\begin{proof}
Зафиксируем точку \(\zeta\in\mathfrak T\), представляющую собой предел
последовательности \(\{f_{z_n}\}_{n\in\mathbb N}\), а также некоторую базисную
окрестность \(U\) этой точки. Введём в рассмотрение три регулярные пары
\(\langle U;\,V_1\rangle\), \(\langle V_1;\,V_2\rangle\) и \(\langle V_2;\,
V_3\rangle\) базисных окрестностей точки \(\zeta\). Наконец, зафиксируем номер
\(k\in\mathbb N\), для которого все члены последовательности
\(\{f_{z_{n+k}}\}_{n\in\mathbb N}\) принадлежат окрестности \(V_3\).

Непосредственно из определений~\ref{prop:2:2:2} и \ref{prop:2:1:1} вытекает
осуществимость алгорифмов \(\mathcal A\) и \(\mathcal B\) со следующими свойствами:
\paragraph\label{prop:2:4:1:1} При любом выборе номера \(n\in\mathbb N\)
применимость алгорифма \(\mathcal A\) к этому номеру означает, что
\(\mathcal A(n)\) есть удовлетворяющий неравенству \(z_n\leqslant\mathcal A(n)\)
элемент множества \(\mathfrak D\), для которого точка \(f_{\mathcal A(n)}\)
не является точкой прикосновения окрестности \(V_2\).
\paragraph\label{prop:2:4:1:2} При любом выборе номера \(n\in\mathbb N\)
значение \(\mathcal B(n)\) определено и представляет собой удовлетворяющий
неравенству \(z_n\leqslant\mathcal B(n)\) элемент множества \(\mathfrak D\).
\paragraph\label{prop:2:4:1:3} Любой номер \(n\in\mathbb N\), для которого
квазиосуществим такой индекс \(\beta\geqslant\mathcal B(n)\), что точка
\(f_{\beta}\) не является точкой прикосновения окрестности \(V_1\), принадлежит
области применимости алгорифма \(\mathcal A\).

\smallskip
Далее будет использоваться стандартная техника, основанная на принципе захвата
(см., напр., \cite[Гл.~9, \S~2.3]{Ku:1973}). Зафиксируем произвольный
арифметический алгорифм \(\mathcal S\) с разрешимым графиком и неразрешимой
областью применимости. Определим действующий на множестве \(\mathbb N^2\) алгорифм
\(\mathcal K\) равенствами
\begin{equation}\label{eq:2:3:2}
	\mathcal K(m,n)\rightleftharpoons\left\{\begin{array}{ll}
		\mathcal A(n)&\text{при } n\eqcirc\mathcal S(m),\\
		z_{n+k}&\text{иначе.}
	\end{array}\right.
\end{equation}
Обозначим через \(\mathfrak M\) множество номеров \(m\in\mathbb N\), для которых
значение \(\mathcal K(m,n)\) определено при всех \(n\in\mathbb N\). Заметим,
что для любого номера \(m\in\mathfrak M\) последовательность \(\{\mathcal K(m,n)
\}_{n\in\mathbb N}\) точек направленного множества \(\mathfrak D\) регулярна
относительно последовательности \(\{z_n\}_{n\in\mathbb N}\) [\eqref{eq:2:3:2},
\ref{prop:2:4:1:1}]. Соответственно, для любого номера \(m\in\mathfrak M\)
или все точки вида \(f_{\mathcal K(m,n)}\) принадлежат окрестности \(V_2\),
или хотя бы одна из таких точек не является точкой прикосновения окрестности
\(V_3\) [\ref{prop:2:2:1}].

Ввиду отмеченных выше особенностей выбора числа \(k\) и построения алгорифма
\(\mathcal A\), сказанное означает возможность указания порождаемого множества
\(\mathfrak N\), содержащего все номера \(m\in\mathbb N\), к которым неприменим
алгорифм \(\mathcal S\), и не содержащего ни одного номера \(m\in\mathbb N\),
для которого определено значение \(\mathcal A(\mathcal S(m))\). При этом, ввиду
неразрешимости области применимости алгорифма \(\mathcal S\), заведомо осуществим
номер \(m\in\mathbb N\), принадлежащий как множеству \(\mathfrak N\), так и области
применимости алгорифма \(\mathcal S\). Вытекающая отсюда неприменимость алгорифма
\(\mathcal A\) к номеру \(\mathcal S(m)\) означает неосуществимость индекса
\(\beta\geqslant\mathcal B(\mathcal S(m))\), для которого точка \(f_{\beta}\)
не была бы точкой прикосновения окрестности \(V_1\) [\ref{prop:2:4:1:3}].
Тем самым, при любом выборе индекса \(\beta\geqslant\mathcal B(\mathcal S(m))\)
точка \(f_{\beta}\) принадлежит окрестности \(U\) [\ref{prop:2:1:1}].

Ввиду произвольности выбора окрестности \(U\), полученные результаты означают
сходимость направленности \(\{f_{\alpha}\}_{\alpha\in\mathfrak D}\) к точке
\(\zeta\).
\end{proof}

\subsection
В качестве простого примера применения утверждения~\ref{cxef:prop} можно
рассмотреть следующую теорему \cite[Теорема~3.2]{Ku:1970} об эквивалентности
слабой и сильной интегрируемости по Риману для всюду на отрезке \([0,1]\)
заданной конструктивной функции:

\subsubsection
{\itshape Пусть функция \(f:[0,1]\to\mathbb R\) такова, что при любом выборе
последовательности \(\{W_n\}_{n\in\mathbb N}\) интегральных дроблений отрезка
\([0,1]\), удовлетворяющей условию\footnote{Здесь через \(\pi(W)\) обозначена
\cite[Введение]{Ku:1970}, \cite[Гл.~7, \S~1]{Ku:1973} измельчённость интегрального
дробления \(W\).}
\begin{equation}\label{eq:2:5:1}
	(\forall n\in\mathbb N)\qquad\pi(W_n)<2^{-n},
\end{equation}
соответствующая последовательность интегральных сумм функции \(f\) окажется
сходящейся. Тогда функция \(f\) интегрируема по Риману.
}

\bigskip
Действительно, легко подбирается такая неубывающая конфинальная последовательность
интегральных дроблений, что всякая регулярная относительно неё последовательность
\(\{W_n\}_{n\in\mathbb N}\) удовлетворяет условию~\eqref{eq:2:5:1}. Необходимая
для возможности приложения утверждения~\ref{cxef:prop} к рассматриваемому случаю
сепарабельность направленности интегральных сумм функции \(f\) тривиальным образом
вытекает из факта непрерывности этой функции \cite[Гл.~5, \S~2,
Теорема~2]{Ku:1973}. Возможность включения произвольной  базисной окрестности
вещественного числа в регулярную пару является следствием теоремы сравнения
\cite[Гл.~2, \S~3, Теорема~20]{Ku:1973}.

\end{document}